\documentclass{amsart}

\usepackage{url}

\usepackage{graphicx}
\usepackage{comment}
\usepackage{amsmath,amssymb,amsthm}
\usepackage{algorithm}
\usepackage{algorithmic}
\usepackage{graphicx,tikz}

\newtheorem*{thm}{Theorem}

\theoremstyle{definition}

\theoremstyle{remark}

\newcommand {\myvec}[1] {{\mbox{\boldmath $#1$}}}

\begin{document}	

	\title{}

\title[]{Refined Least Squares for Support Recovery}
\thanks{S.S. is supported by the NSF (DMS-2123224) and the Alfred P. Sloan Foundation.}

\author[]{Ofir Lindenbaum}
\address{Program in Applied Mathematics, Yale University, New Haven, CT 06511, USA}
\email{ofir.lindenbaum@yale.edu}

\author[]{Stefan Steinerberger}
\address{Department of Mathematics, University of Washington, Seattle, WA 98195, USA}
\email{steinerb@uw.edu}
	
	\begin{abstract}
		We study the problem of exact support recovery based on noisy observations and present Refined Least Squares (RLS). Given a set of noisy measurement 
$$ \myvec{y} = \myvec{X}\myvec{\theta}^* + \myvec{\omega},$$
and $\myvec{X} \in \mathbb{R}^{N \times D}$ which is a (known) Gaussian matrix and $\myvec{\omega} \in \mathbb{R}^N$ is an (unknown) Gaussian noise vector, our goal is to recover the support of the (unknown) sparse vector $\myvec{\theta}^* \in \left\{-1,0,1\right\}^D$. To recover the support of the $\myvec{\theta}^*$ we use an average of multiple least squares solutions, each computed based on a subset of the full set of equations. The support is estimated by identifying the most significant coefficients of the average least squares solution. We demonstrate that in a wide variety of settings our method outperforms state-of-the-art support recovery algorithms. 
	\end{abstract}

	\maketitle

%	\IEEEpeerreviewmaketitle

	\section{Introduction}

	\subsection{The Problem}
The task of estimating the support of a sparse signal based on noisy measurements plays a key role in many applications, such as, image denoising \cite{elad2006image,zhao2014hyperspectral}, communication \cite{qin2018sparse}, machine learning \cite{guo2010action,yang2011learning,DUFS}. Here, we consider support recovery of an unknown ternary signal $\myvec{\theta}^* \in \left\{-1,0,1\right\}^D$, given the noisy vector
	$$ \myvec{y} = \myvec{X}\myvec{\theta}^* + \myvec{\omega},$$ 
	where $\myvec{X} \in \mathbb{R}^{N \times D}$ is a measurement matrix whose entries are independent random variables with mean $0$ and variance $1$, and the noise $\myvec{\omega} \in \mathbb{R}^N$ is Gaussian. Exploiting the known sparsity $k$ of the unknown signal $\myvec{\theta}^*$, our goal is to recover its support, i.e. all indices $d$ such that $\myvec{\theta}^*_d \neq 0$. We are interested in the case when the system is underdetermined and has more variables than equations $D>N$.

	\begin{center}
\begin{figure}[h!]
\begin{tikzpicture}[scale=0.5]
\node at (-7,0) {$y=$};
\node at (0,0) {\includegraphics[scale=0.5]{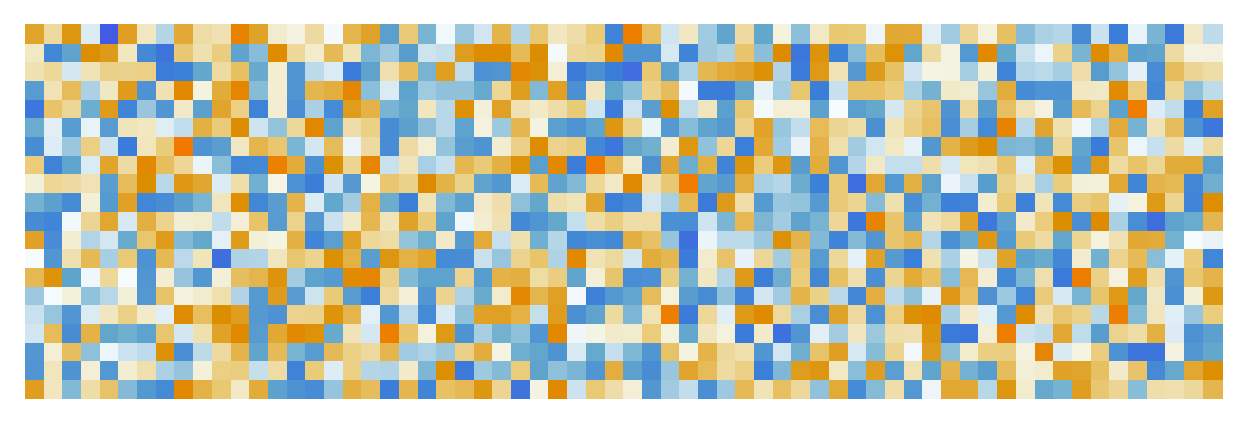}};
\node at (7.1,0) {$\theta + w$};
\end{tikzpicture}
\caption{We try to recover the support of $\theta$ from the observations $X$ and $y$, where $y=X\theta + \omega$. The (known) matrix $X$ is a Gaussian random matrix, so is the (unknown) noise $\omega$, we try to recover the support of $\theta$ with few measurements.}
\end{figure}
\end{center}
	
\subsection{Existing Results.}	A popular approach for recovering the support of $\myvec{\theta}^*$ is using an $\ell^1$ regularized quadratic programming problem, also known as the Least Absolute Shrinkage and Selection Operator (LASSO) \cite{lasso}. The LASSO is cast as a convex problem, and can be solved using efficient optimization schemes \cite{lasso2, lasso3}. Several authors have introduced iterative algorithms improving the support recovery capabilities of the LASSO, these include: Iterative Re-weighted Least Squares (IRLS) \cite{irls}, Iteratively Re-weighted $\ell^1$ minimization (IRL1) \cite{irl1}, and Iterative Support Detection (ISD) \cite{isd}. Other solutions for support recovery include greedy algorithms such as Orthogonal Matching Pursuit (OMP) \cite{omp}, and its extensions \cite{elad2009plurality,stomp,cosamp}, or non convex schemes such as the Trimmed LASSO (TL) \cite{tl}, smoothly clipped absolute deviation (SCAD) \cite{scad}, or stochastic gates (STG) \cite{stg}. For a more complete survey of the existing methodologies, we refer the reader to \cite{arjoune2017compressive,bruckstein2009sparse,marques2018review,mousavi2019survey}.
	Following the approach of Randomly Aggregated Least Squares (RAWLS) proposed by the authors in \cite{Rawls}, we introduce Refined Least Squares (RLS), a procedure for estimating the support of $\myvec{\theta}^*$ based on the noisy vector $\myvec{y}$ and the observations matrix $\myvec{X}$. Our algorithm is greedy, relying on estimating the support (and sign) of one coefficient at a time. We demonstrate that Refined Least Squares (RLS) effectively recovers the support, specifically in the low information regime where $N$ is much smaller than $D$, the regime where the noise is high and there is a large number of active coefficients ($k$).
	
	\section{Refined Least Squares}
	\subsection{RAWLS}
	The purpose of this section is to introduce the main idea behind Refined Least Squares (RLS). Our starting point is RAWLS (\textbf{R}andomly \textbf{A}ggregated un\textbf{W}eighted \textbf{L}east \textbf{S}quares), a method introduced by the authors in \cite{Rawls}. The main idea behind RAWLS is somewhat counter-intuitive: we are given $N$ equations with $D$ variables (with special emphasis on $N < D$)
	$$ \myvec{y} = \myvec{X}\myvec{\theta}^* + \myvec{\omega}.$$
	It is easy to see that 
	$$ \myvec{\widehat{\theta}} = \arg\min_{\theta} \| \myvec{X} \myvec{\theta} - \myvec{y} \|$$ 
	does not usually lead to very good results (in particular, $\myvec{\widehat{\theta}}$ is not typically sparse). The main idea behind RAWLS is to do a least squares regressions on a \textit{subset} of the equations: for any arbitrary subset $A \subset \left\{1,2,\dots, N\right\}$, we use $\myvec{X}_A$ to denote the reduction of $\myvec{X}$ to those rows indexed by $A$ and likewise for the vector $\myvec{y}_A$. We can then consider
	$$ \myvec{\widehat{\theta}}_A = \arg\min_{\theta \in \mathbb{R}^D} \| \myvec{X}_A \myvec{\theta} - \myvec{y}_A \|.$$ 
	This will typically lead to even worse results compared to using least squares on the full system since we work with less information. However, since there are \textit{many} different subsets $A \subset \left\{1,2,\dots, N\right\}$, we can obtain several different estimators $ \myvec{\widehat{\theta}}_{A_1},   \myvec{\widehat{\theta}}_{A_2}, \dots \myvec{\widehat{\theta}}_{A_m}$ and average them. RAWLS chooses $m$ random subsets of size $|A_i| = \left\lfloor0.6\min(D,N) \right\rfloor$ and averages over the different least squares solutions. The main observation that makes RAWLS work is that the coordinate (and the sign) of these averages with the largest (absolute) entry likely corresponds to a coordinate for which $\myvec{\theta}^*$ has a nonzero entry (with the same sign). This coordinate is removed in a greedy fashion and the procedure is reiterated on the remaining problem. 
	
	\subsection{Refined Least Squares} Refined Least Squares (RLS) is a significant refinement over RAWLS and leads to a very effective algorithm (see \S 3 for results). As was already pointed out in \cite{Rawls}, we believe that the idea of suitably averaging over least square solutions computed over \textit{subsets} of equations to be a promising idea and not yet suitably explored. We present two additional ideas which we call 
	\begin{enumerate}
	\item \textit{averaged guessing} and 
	\item the \textit{OMP flip}
	\end{enumerate}
	They lead to Refined Least Squares (RLS) which achieves state-of-the-art results in several challenging regimes. We now explain these two new ideas.
	
	\subsubsection{Averaged Guessing.} RAWLS produces a guess for a coordinate with nonzero support by averaging least square recoveries over random subsets of size  $|A_i| = \left\lfloor0.6\min(D,N) \right\rfloor$. As demonstrated in Fig. \ref{fig:m-eval}, the constant $0.6$ leads to the decent results in terms of support recovery but not in a very strict sense -- any constant in the range $[0.6, 0.9]$ seems to lead to rather comparable behavior (this is related to the distribution of singular values of random matrices, see below). The main new idea is to run several such predictions, for $|A_i| \in [ 0.85\min(D_k,N), 0.9\min(D_k,N)]$, where $D_k$ is the number of variables at peeling step $k=1,...,D$. Then, using a majority vote among the top selected components, we get an additional layer of stability leading to better results. 
	
	\subsubsection{The OMP flip.} The second new ingredient is based on the following idea. As RAWLS proceeds in greedily selecting coordinates, the problem tends to become more difficult as the leading variables (corresponding to the strongest coefficients) are removed because the signal to noise ratio shrinks. In particular, after $k-1$ steps, we have a system of $N$ equations for $D - k +1$ unknowns which is mainly dominated by noise once $k \sim D$. The final step is as follows: for $D$ Gaussian vectors $\myvec{g}_1, \dots, \myvec{g}_D \in \mathbb{R}^N$ (the columns of the matrix $\myvec{X}$), we are given $\myvec{y} = (\pm 1) \cdot \myvec{g}_i + \myvec{\omega}$ and are asked to reconstruct $i$ from knowing $\myvec{y}$ and the $D$ vectors $\myvec{g}_1, \dots, \myvec{g}_D \in \mathbb{R}^N$. It is difficult to imagine an estimator better than
	$$ \widehat{i} = \arg \max_{1 \leq i \leq D} \left| \left\langle \myvec{y}, \myvec{g}_i \right\rangle \right|.$$
	This way of reconstructing coordinates is the main idea behind Orthogonal Matching Pursuit (OMP) \cite{omp}. Once there is a single coordinate left to reconstruct, (RLS) switches to this OMP heuristic. This also leads to additional stability. 
	
	\textit{The Algorithm.}
	To estimate the support we use an ensemble of estimates, each based on a different subset size $|A_i|=n_0$. We fix the subset size 
	$$n_0 \in [ 0.85\min(D_k,N), 0.9\min(D_k,N)]$$
	 and start the peeling procedure to estimate the set of support indices ${\mathcal I}\subset \{1,...,D\}$. Using $m$ (we use $m=100$) random subsets of the equations (of size $n_0$ each), the RLS estimate is based on the following approximation 
	 $$ \myvec{\bar{\theta}} = \frac{1}{m} \sum_{i=1}^{m}\myvec{\widehat{\theta}}_{A_i}.$$ We estimate a support coefficient as 
	$\widehat{\theta}_{\ell}= \text{sign}(\max{\myvec{\bar{\theta}}}),$ where $\ell=\arg\max{|{\myvec{\bar{\theta}}}|} $ is the index of the support variable. The process is repeated $5$ times with different subset size $n_0 \in [ 0.85\min(D_k,N), 0.9\min(D_k,N)].$
	Finally, the set of support indices $\mathcal I$ is estimated as using a majority voting over the different sets identified in the $5$ peeling procedures: we pick the coordinate most frequently identified as most likely (in case of a tie, we choose randomly among the most frequently occuring ones).  After that, we remove $\myvec{X}^{\ell}$ from $\myvec{X}$, where $\myvec{X}^{\ell}$ is the $\ell$'s column of the matrix, and update $\myvec{y}$ by subtracting $\myvec{X}^{\ell}\widehat{\theta}_{\ell}$. This process is repeated $k-1$ times for estimating the leading support indices, the last coefficient is estimated using OMP. 
	A description of the peeling procedure is also presented in Algorithm \ref{alg:pseudocode}.

	\begin{algorithm}[h]
		\caption{Refined Least Squares (RLS)}
		\label{alg:pseudocode}
		
		{\bfseries Input:} Observations matrix $\myvec{X}$, and target vector $\myvec{y}$.\\
		%	{\bfseries Output:} Support set ${\cal I}\subset \{1,...,D\}$.
		\begin{algorithmic}[1]
			\FOR {$d=1:k-1$}
			%	\FOR{ $S\in {0.7D,0.75D,0.8D}$}
			\FOR{$i=1:m$}
			\STATE Draw a set of equations $A_i \subset \{1,...,N\}$ of size $n_0$.
			\STATE Use least squares to find $\bar{\myvec{\theta}}_{A_i}=\myvec{X}_{A_i}^{\dagger}\myvec{y}_{A_i}$ 
			\ENDFOR
			\STATE Compute the approximation $ \myvec{\bar{\theta}} = \frac{1}{m} \sum_{i=1}^{m}\myvec{\widehat{\theta}}_{A_i}$
			\STATE Estimate $\widehat{\theta}_{\ell}= \text{sign}(\max{\myvec{\bar{\theta}}})$, where $\ell=\arg\max{|{\myvec{\bar{\theta}}}|} $
			\STATE Remove the $\ell$ column from $\myvec{X}$ and update $\myvec{y}$ by subtracting $\myvec{X}^{\ell}\widehat{\theta}_{\ell}$, where $\myvec{X}^{\ell}$ is the $\ell$'s column of the matrix.
			\ENDFOR
			%	\ENDFOR
			\STATE Estimate the $k$'th coefficient using OMP.
		\end{algorithmic}
	\end{algorithm}
	%\textbf{[HOW BIG IS $M$ IN PRACTICE?]}\\
	
	\subsection{The Theorem.}
	Our main theoretical contribution is a rigorous result explaining why, for this type of problem, it is actually advantageous to use a least squares approach over a \textit{reduced} set of equations (something that is maybe counter-intuitive since we `throw away' information). We will now state the theorem which is relevant for a large number of algorithms of this type. Our setting is as follows: we assume that $N < D$ is given and that $\myvec{X} \in \mathbb{R}^{N \times D}$ is a random Gaussian matrix with each entry being an independent $\mathcal{N}(0,1)$ random variable. Let us assume that the vector $\myvec{\theta}^* \in \mathbb{R}^D$ is fixed and that $\myvec{\omega} \in \mathbb{R}^D$ is a Gaussian perturbation. We are given $\myvec{y} = \myvec{X \theta^* + \omega}$. We compute a least squares approximation $\myvec{X}^{\dagger} \myvec{y}$, where $\myvec{X}^{\dagger}=\myvec{X}^T(\myvec{X X}^T)^{-1}$, and are interested in how close this approximation is to $\myvec{X}^{\dagger} \myvec{\theta}^*$. This means that the main question is regarding the size of $ \|\myvec{X}^{\dagger} \myvec{y} - \myvec{X}^{\dagger} \myvec{X} \myvec{\theta}^* \|$: we hope that the least squares reconstruction of $\myvec{y}$ is very close to that of $ \myvec{X} \myvec{\theta}^*$. The subsequent Theorem tells us that how close these two numbers are depends on the dimensions of $\myvec{X} \in \mathbb{R}^{N \times D}$: their distance grows as $N/D$ gets closer to 1.
	
	\begin{thm}
		If the size of the matrix $\myvec{X}$ tends to infinity, with $N/D < 1$ fixed, and $\omega_i \sim \mathcal{N}(0,1)$, we have
		$$ \mathbb{E}_{{\myvec{X}, \myvec{\omega}}}~ \| \myvec{X}^{\dagger} \myvec{y} - \myvec{X}^{\dagger} \myvec{X} \myvec{\theta}^* \| = (1+ o(1)) \sqrt{\frac{N}{D-N}}.$$
		
	\end{thm}

\begin{figure}[htb!]
		
		\centering
		\includegraphics[width=0.8\textwidth]{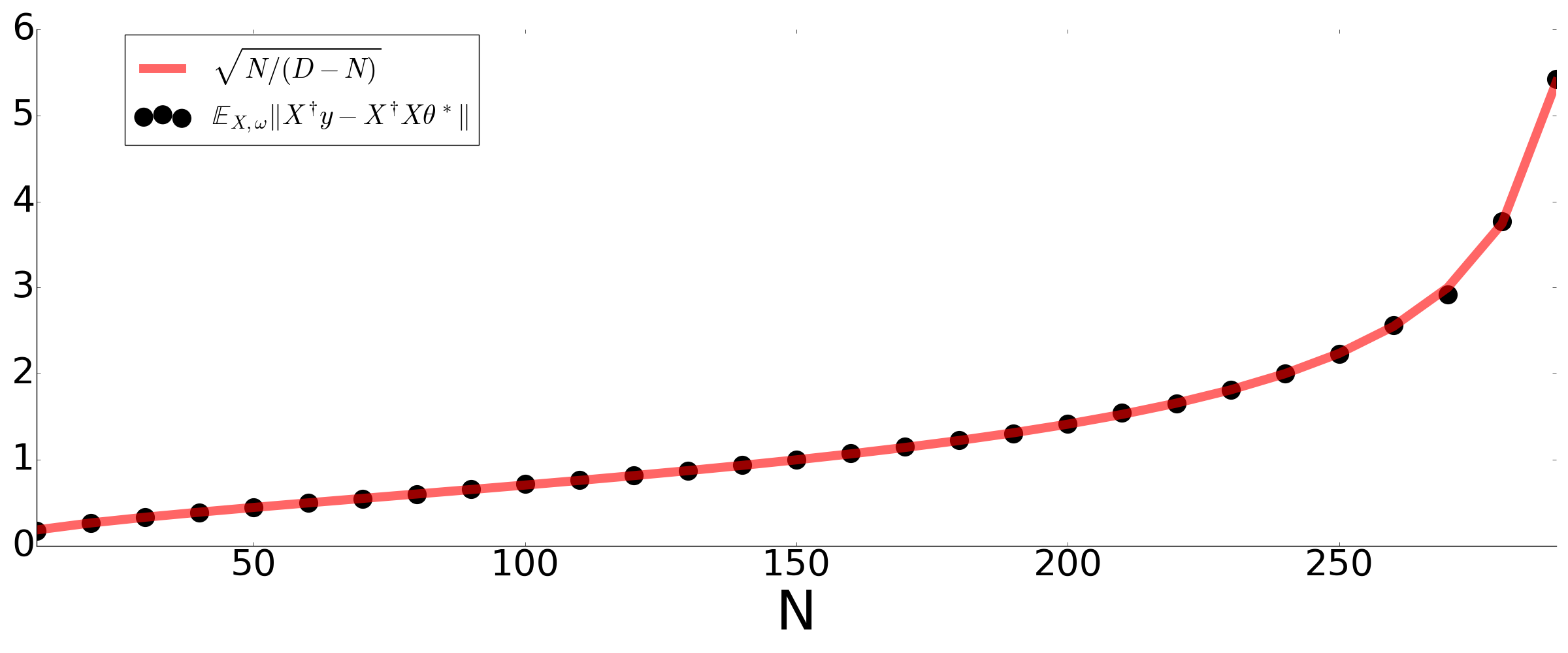}
		
		\vskip -0.1  in
		\caption{Numerical evaluation of the relation predicted by Theorem 1 for $D=300$ and $10 \leq N \leq 290$. As the number of equations $N$ tends to $D$ (the number of variables), the expected $\ell^2$ norm of $\myvec{X}^{\dagger}\myvec{y}- \myvec{X}^{\dagger}\myvec{X}\myvec{\theta}^*$ (black dots) grows like $N^{1/2} (D-N)^{-1/2}$ (red line). }
\vskip -0.05  in		
		\label{fig:error}
	\end{figure}

	We note that $\mathbb{E} \| \myvec{X}^{\dagger} \myvec{X}\myvec{\theta}^* \|$ may be thought of as the projection onto a random hyperplane and we may thus assume that
	$$ \mathbb{E} \|\myvec{X}^{\dagger}\myvec{X}\myvec{\theta}^* \| \sim \sqrt{\frac{N}{D}} \|\myvec{X} \myvec{\theta}^*\|.$$
	This shows that both $\| \myvec{X}^{\dagger}\myvec{X}\myvec{\theta}^* \|$ and $ \| \myvec{X}^{\dagger} \myvec{y} - \myvec{X}^{\dagger}\myvec{X} \myvec{\theta}^* \|$ have a priori same dependence on the ratio of the dimensions of $\myvec{X}$. $N/D$ should not be too small because this amounts to throwing away a lot of information. However, it also should not be too close to 1 either since in that case the size of the error starts growing dramatically: this is because the matrix $\myvec{X}$ suddenly starts having singular values closer to $0$ which badly influence the inversion problem. The main ingredient is the Marchenko-Pastur Theorem which holds for general random rectangular matrices whose random variables are independent and have mean value $0$ and variance $1$:  the Theorem also applies to measurement matrices that are not Gaussian (as we demonstrate in Section \ref{sec:exp}). A numerical evaluation supporting the Theorem appears in Fig. \ref{fig:error}.

	\begin{figure}[htb!]
		
		\centering
				\includegraphics[width=0.8\textwidth]{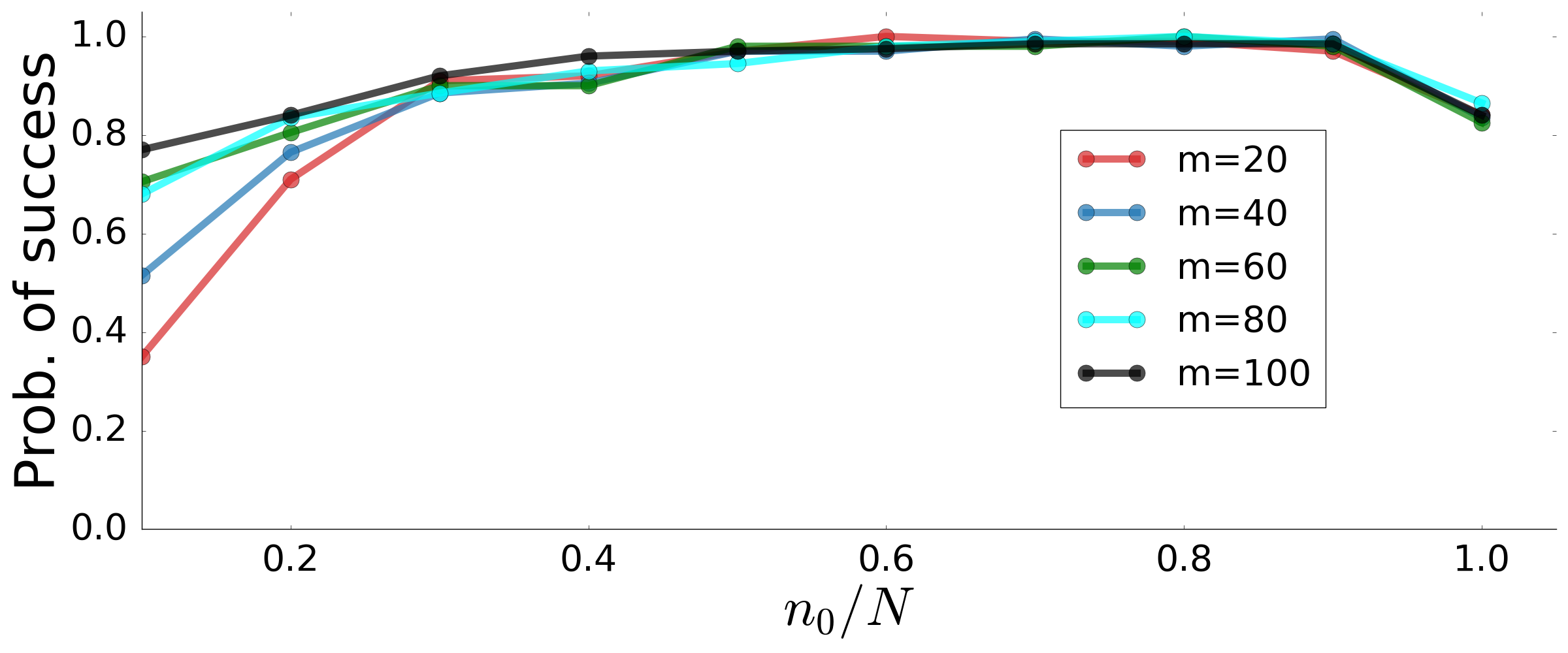}

		\caption{Probability of successfully recovering the support of the signal using RLS with a subset of size $n_0$ (x-axis presents the normalized subset size $n_0/N$). In this experiment, we use $D=64$ variables $N=50$ equations and a sparsity of $k=20$. We use several values for $m$, the number of subsets used for estimating the support via Algorithm 1. }
		
		\label{fig:m-eval}
	\end{figure}

	\section{Experimental Results}\label{sec:exp}
	In the following subsection we evaluate the support recovery capabilities of RLS in different settings. We compare the performance of RLS to several strong baselines, such as: the LASSO \cite{lasso}, IRL1 \cite{irl1}, TL \cite{tl}, OMP \cite{omp}, ISD \cite{isd} and RAWLS \cite{Rawls}. We apply RLS with $m=100$ which was shown stable across different settings (as demonstrated in Fig. \ref{fig:m-eval} smaller values of $m$ may also work in practice). We evaluate performance in terms of the empirical probability for perfect support recovery. This probability is estimated as the portion of simulations which obtained perfect support recovery out of $200$ runs. Perfect support recovery is counted only if $S(\myvec{\theta})=S(\widehat{\myvec{\theta}})$, where $S(\myvec{\theta}):= \{i \in 1,...,D | \myvec{\theta}_i \neq 0 \}$.

\begin{figure}[htb!]	
		\centering
		\includegraphics[width=0.8\textwidth]{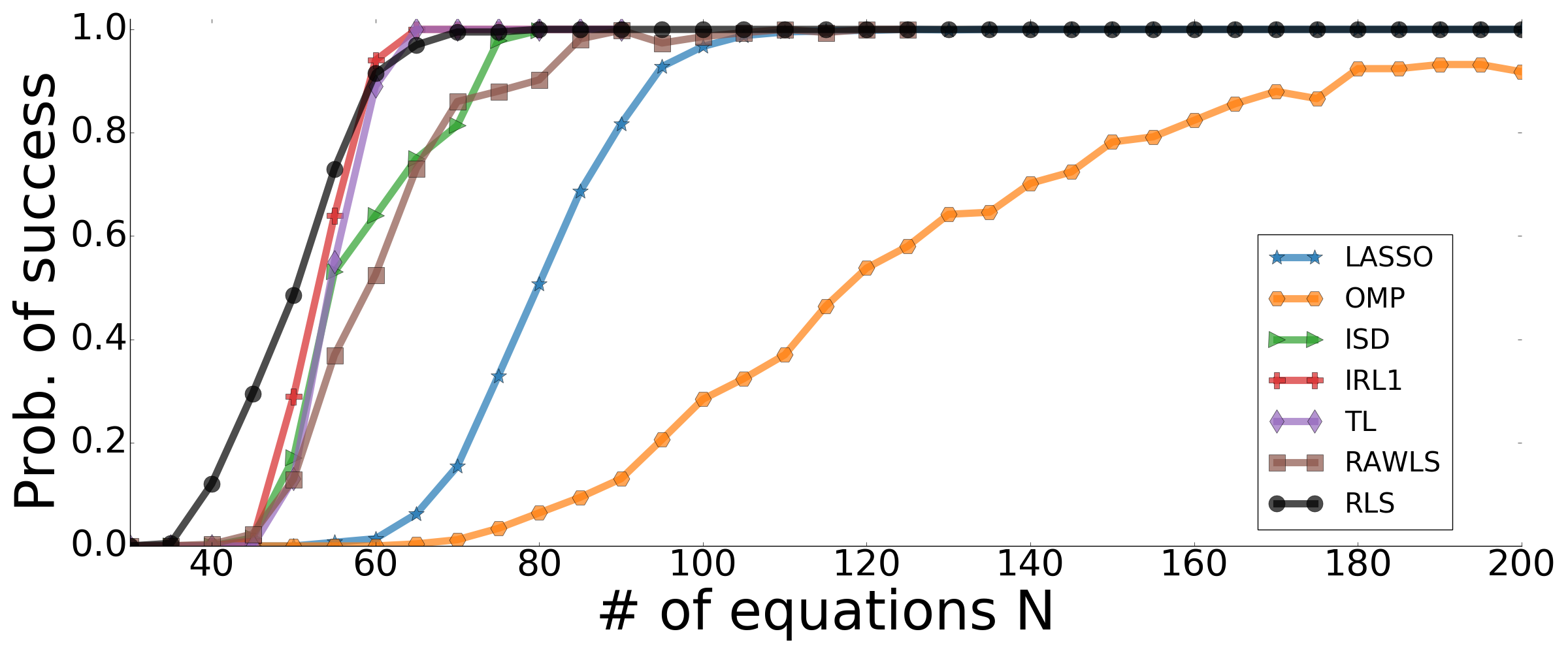}
		\includegraphics[width=0.8\textwidth]{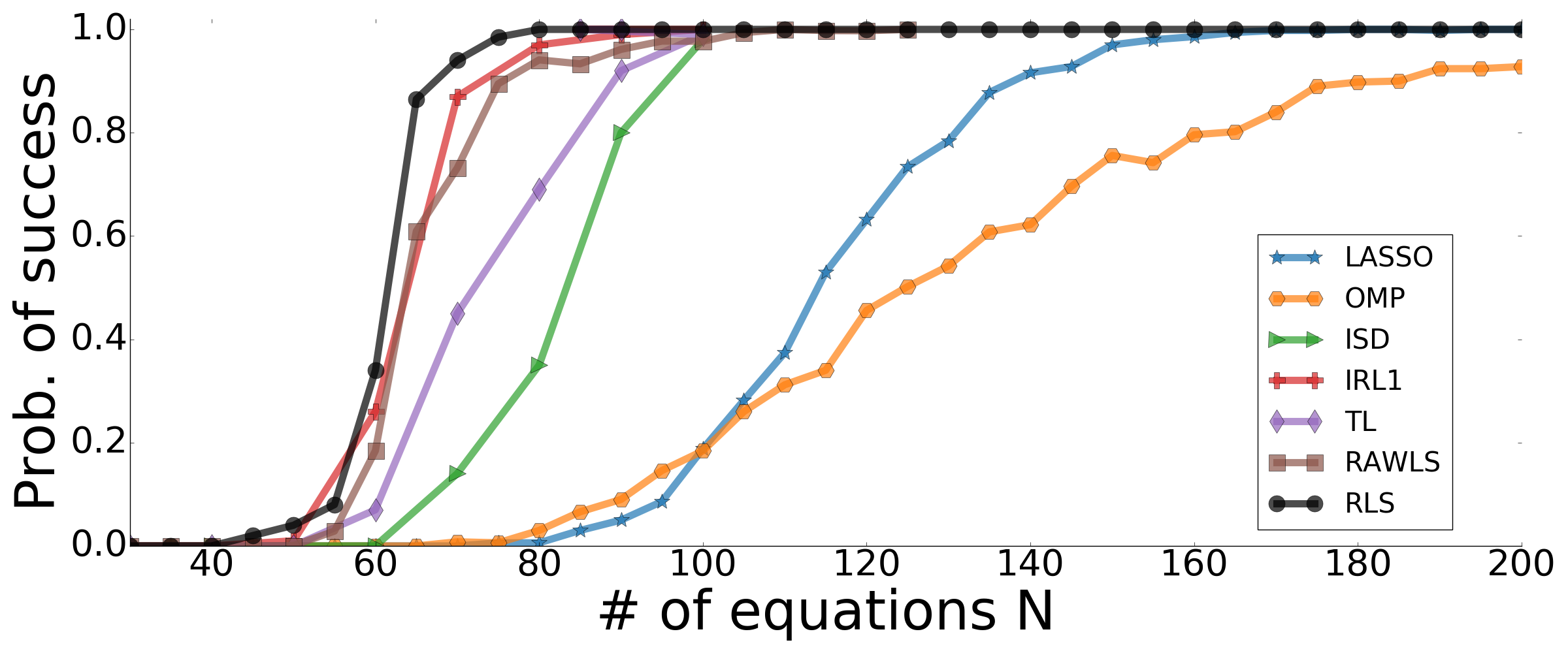}	
		\caption{Evaluating the probability of exact support recovery vs. the number of measurements $N$. We use $D=64$ variables, with a sparsity $k=30$, and a design matrix $\myvec{X}$ with values drawn from $N(0,1)$. We compare the proposed procedure (RAWLS) to several baselines for:  $\sigma=0.5$ (top panel) and $\sigma=1$ (bottom panel).}
		
		\label{fig:k30}
	\end{figure}

\begin{figure}[htb!]
		
		\centering
		\includegraphics[width=0.8\textwidth]{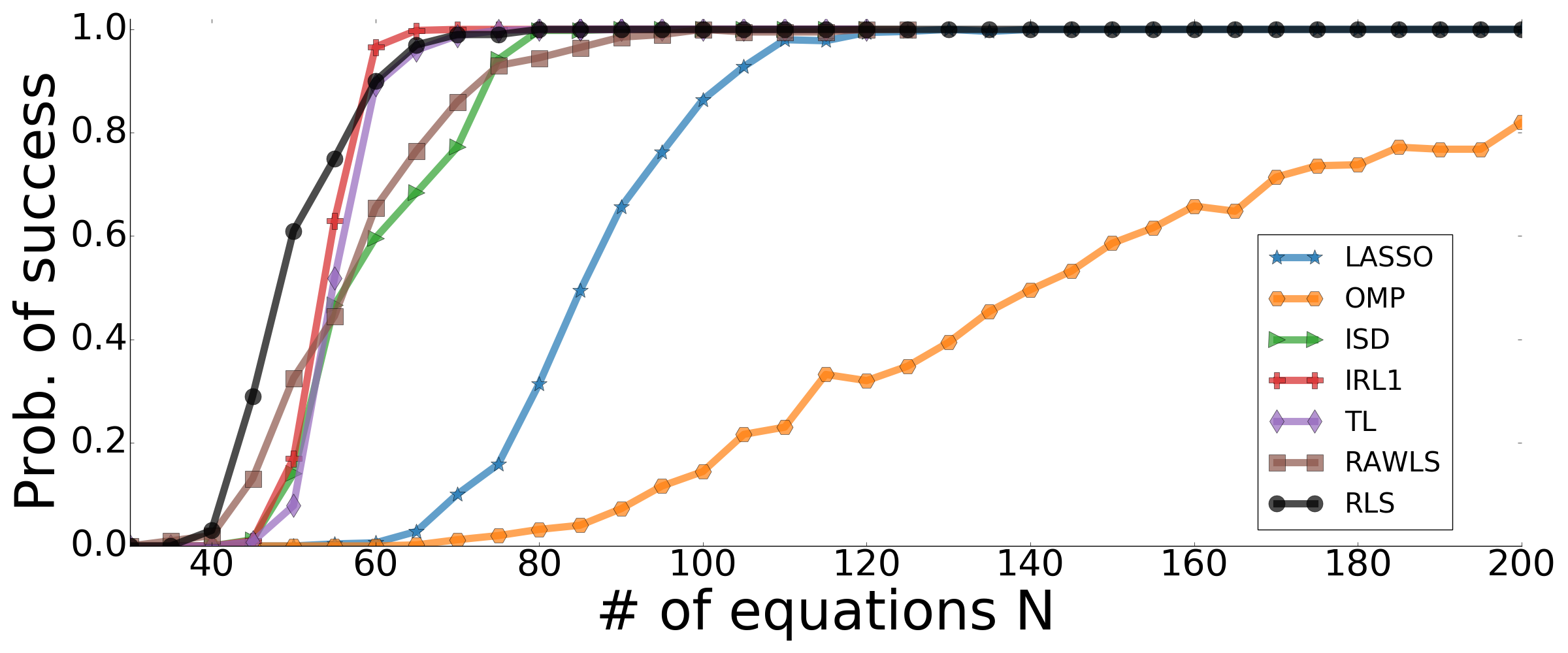}
		\includegraphics[width=0.8\textwidth]{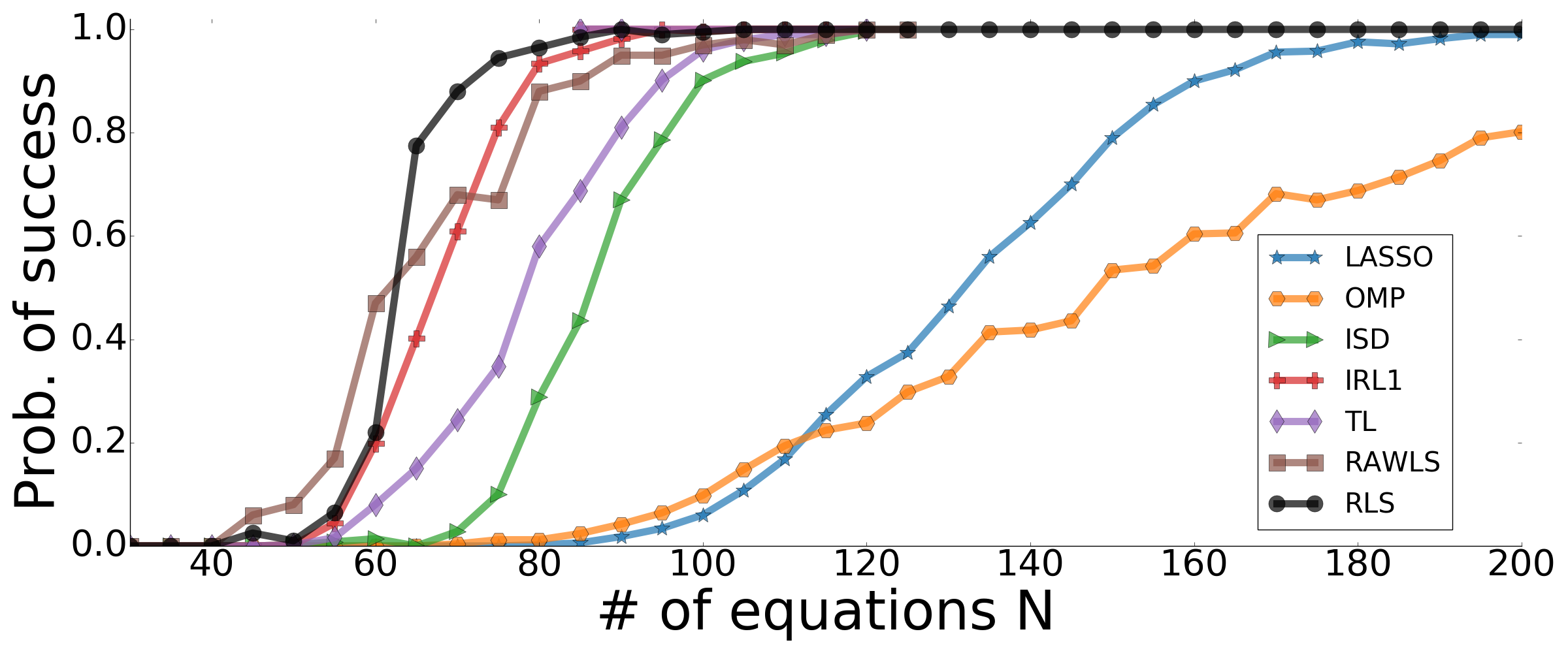}		
		\caption{Evaluating the probability of exact support recovery vs. the number of measurements $N$. We use $D=64$ variables, with a sparsity $k=30$, and a design matrix $\myvec{X}$ with values drawn from $N(0,\myvec{\Sigma})$, where $\Sigma_{i,j}=0.3^{|i-j|}$. We compare the proposed procedure (RAWLS) to several baselines for:  $\sigma=0.5$ (top panel) and $\sigma=1$ (bottom panel).}	
		\label{fig:k30-corr}
	\end{figure}

\begin{figure}[htb!]		
		\centering
				\includegraphics[width=0.8\textwidth]{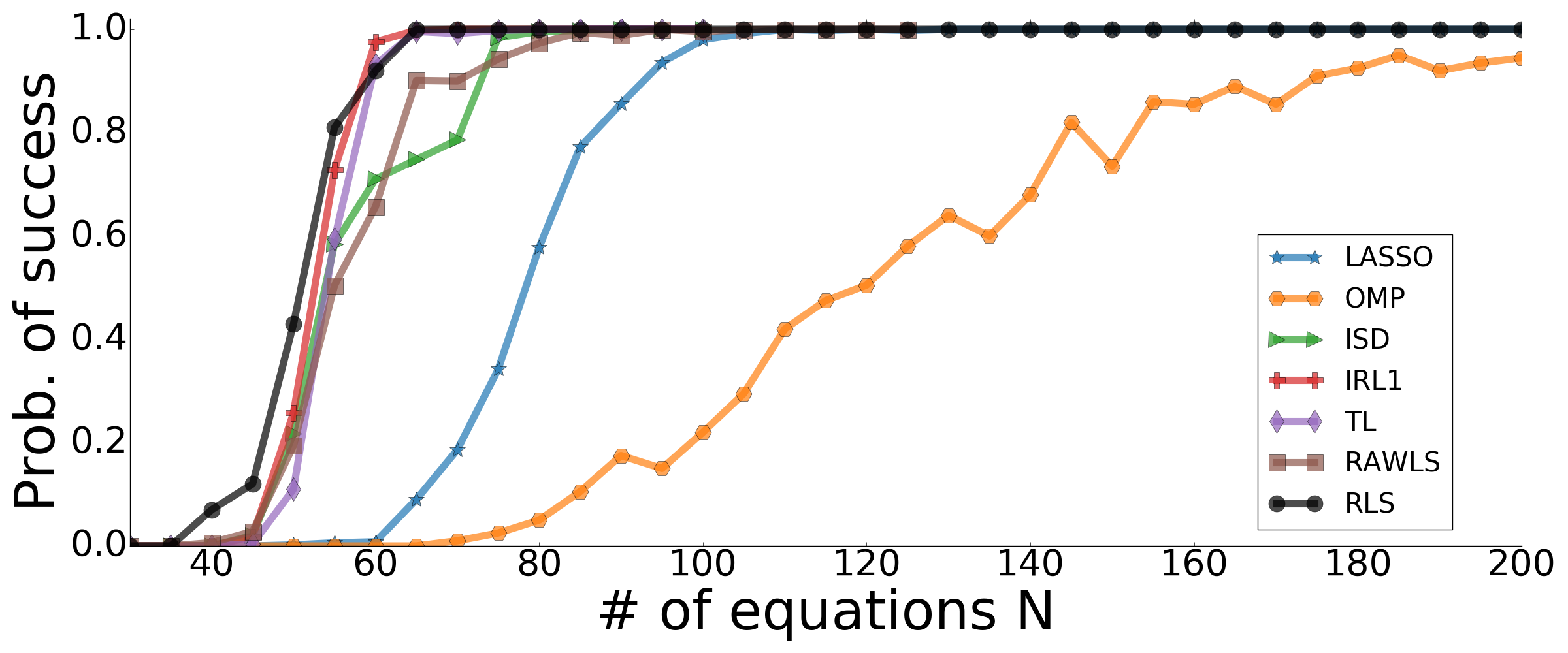}
		\includegraphics[width=0.8\textwidth]{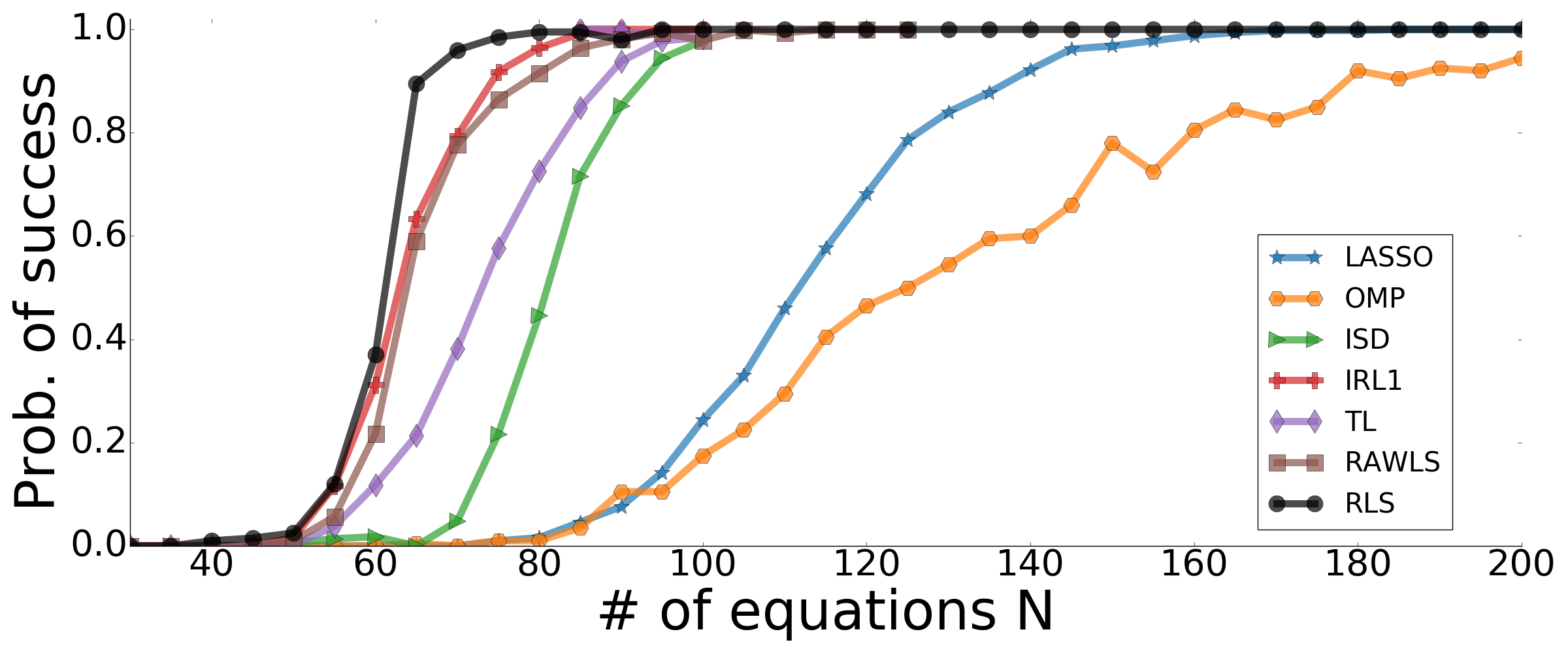}	
		\caption{Evaluating the probability of exact support recovery vs. the number of measurements $N$. We use $D=64$ variables, with a sparsity $k=30$, and a design matrix $\myvec{X}$ with values drawn from a fair Bernoulli distribution. We compare the proposed procedure (RAWLS) to several baselines for:  $\sigma=0.5$ (top panel) and $\sigma=1$ (bottom panel).}
		
		\label{fig:k30-binary}
	\end{figure}

\begin{figure}[htb!]
		
		\centering
		\includegraphics[width=0.8\textwidth]{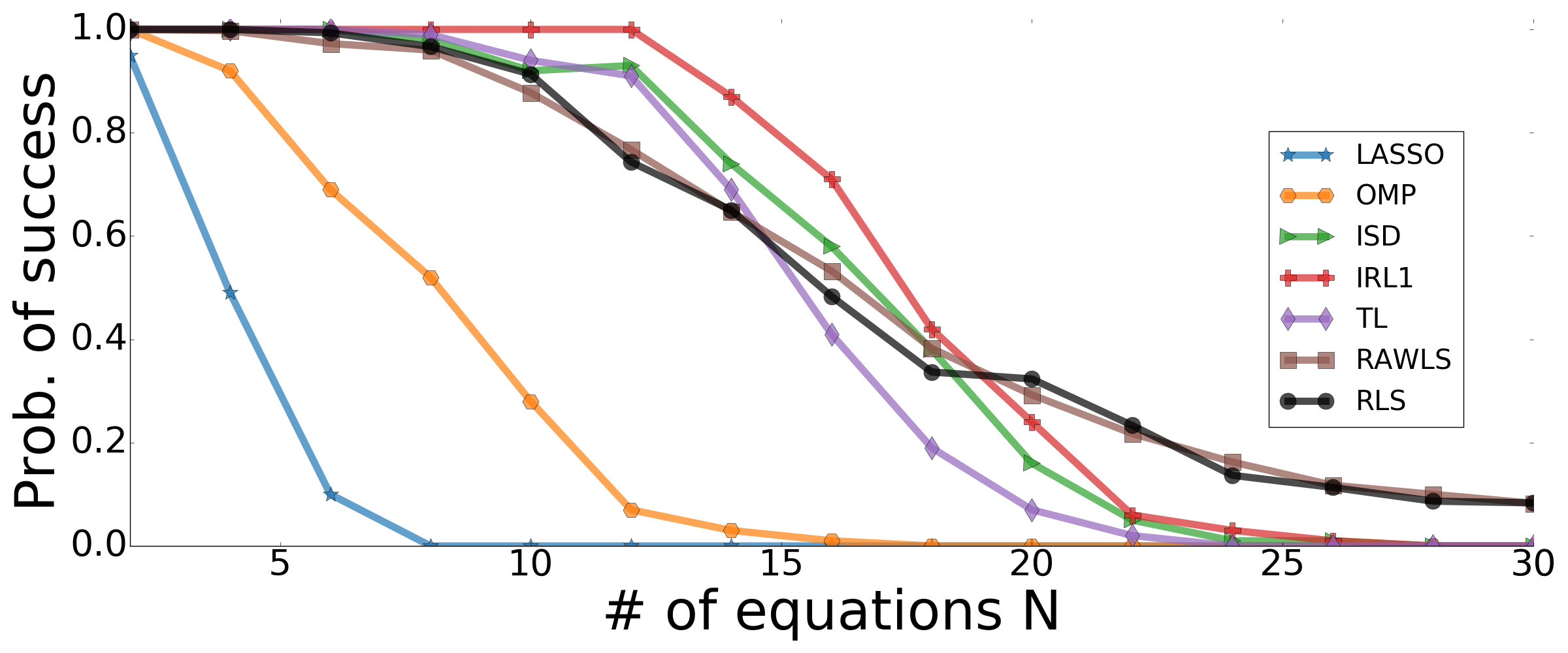}

		\caption{Evaluating the probability of exact support recovery vs. sparsity level $k$. We use $D=64$ variables, with $N=40$ observations, $\sigma=1$ and different sparsity levels.}
		
		\label{fig:sparsity}
	\end{figure}

	For the first example (see Fig. \ref{fig:k30}), we use a design matrix $\myvec{X}$ with values drawn independently from $N(0,1)$, the number of variables is $D=64$ with sparsity $k=30$. This example demonstrates the advantage of RLS over several strong baselines for different noise levels ($\sigma=0.5$ and $\sigma=1$).	
	Then (see Fig. \ref{fig:k30-corr}), we repeat the experiment but use a correlated design matrix $\myvec{X}$. To generate $\myvec{X}$, we first construct a Toeplitz covariance matrix $\myvec{\Sigma}$ with values $\Sigma_{i,j}=0.3^{|i-j|}$, where $i,j=1,...,D$, then we generate $\myvec{X}$ by drawing its values from $N(0,\myvec{\Sigma})$. 	
	In the next experiment (see Fig. \ref{fig:k30-binary}), we investigate what happens when the design matrix $\myvec{X}$ takes its entries from a fair Bernoulli distribution. Each entry $\myvec{X}_{ij}$ is taken independently and uniformly at random from $\left\{-1,1\right\}$.	
	In the last experiment (see Fig. \ref{fig:sparsity}), we investigate the behavior of these methods with respect to sparsity. We fix $D=64$ variables, $N=40$ equations and assume the error on the right-hand side to be of scale $\sigma = 1$. We observe that RLS can recover the true support with some nonzero (albeit small) likelihood even when the solution is far from sparse ($k=30$).
	
	\section{Proof of the Theorem}

	\begin{proof}
		We will use $\myvec{X}^{\dagger} \myvec{y}$ to denote the $\ell^2-$smallest vector satisfying
		$$ \myvec{X}^{\dagger} \myvec{y}= \arg \min_{{\small \myvec{\theta}} \in \mathbb{R}^D} \| \myvec{X \theta - y}\|.$$
		This solutions is contained in the vector space $ V$ spanned by the rows of the matrix (if $\myvec{\theta}$ had a component that was orthogonal to these rows, then it would not have any effect in the matrix multiplication $\myvec{X} \myvec{\theta}$ and removing that component would result in a smaller $\ell^2-$norm). Since the number of variables, $D$, is larger than the number of equations, $N$, and $\myvec{X}$ is Gaussian we know that the minimum is $0$ with likelihood $1$ (the system is underdetermined and thus has a solution with likelihood $1$). By linearity, we have  $\myvec{X}^{\dagger} (\myvec{y} - \myvec{X} \myvec{\theta}^* ) = \myvec{X}^{\dagger} \myvec{\omega}$
		and thus
		$ \| \myvec{X}^{\dagger} \myvec{y} -  \myvec{X}^{\dagger} \myvec{X} \myvec{\theta}^* \| = \| \myvec{X}^{\dagger}  \myvec{\omega}^*\|.$
		At this point, we fix the random matrix $\myvec{X}$ and take an expectation over $\myvec{\omega}$. We observe that each entry of $\myvec{\omega}$ is an independent $\mathcal{N}(0,1)$ Gaussian. Let us now assume that the singular value decomposition of $\myvec{X}$ is given by
		$ \myvec{X} = \myvec{U} \myvec{\Sigma V}^*$ and thus $\myvec{X}^{\dagger} = \myvec{V \Sigma}^{\dagger} \myvec{U}^*.$
		The matrices $\myvec{U}$ and $\myvec{V}$ are unitary and the distribution of Gaussian vectors is invariant under an application of orthogonal matrices, thus
		$ \mathbb{E}_{{\omega}} \| \myvec{X}^{\dagger}  \myvec{\omega}\| =  \mathbb{E}_{{\omega}} \| \myvec{\Sigma}^{\dagger}  \myvec{\omega}\|.$
		Writing $\myvec{\omega} = \left( \omega_1, \dots, \omega_{N}\right)$ and writing $\sigma_1, \dots, \sigma_{N}$ for the singular values of $\myvec{X}$, we can use the explicit form of $\myvec{\Sigma}^{\dagger} $ to compute
		$$  \mathbb{E}_{{\omega}} \| \myvec{\Sigma}^{\dagger}  \myvec{\omega}\| =  \mathbb{E}_{{\omega}} \left( \sum_{i=1}^{N} \frac{ \omega_i^2}{\sigma_i^2} \right)^{1/2}.$$
	Jensen's inequality leads to
		\begin{align*}
			\mathbb{E}_{{\omega}} \left( \sum_{i=1}^{N} \frac{ \omega_i^2}{\sigma_i^2}  \right)^{1/2} &\leq \left(\mathbb{E}_{{\omega}}  \sum_{i=1}^{N} \frac{ \omega_i^2}{\sigma_i^2}  \right)^{1/2} \\
			&=  \left( \frac{1}{N} \sum_{i=1}^{N} \frac{ 1}{\sigma_i^2} \right)^{1/2} \sqrt{N}.
		\end{align*}
		We note that this inequality is strict, however, standard deviation inequalities imply that the inequality becomes asymptotically sharp because the sum starts to tightly concentrate around its mean value as the dimension of the matrix increases.
		The distribution of singular values of rectangular matrices of dimension $N \times D$ is given by the Marchenko-Pastur distribution with parameter $0 \leq N/D \leq 1$. Using $f_{N/D}$ to denote the density of the Marchenko-Pastur distribution, we see that, in the limit as $N, D \rightarrow \infty$ with $N/D$ fixed,
		$$\mathbb{E}_{{X}}  \frac{1}{N} \sum_{i=1}^{N} \frac{ 1}{\sigma_i^2} = (1+ o(1)) \frac{1}{D}  \int_{0}^{\infty} \frac{ f_{N/D}(x)}{x} dx.$$

		It remains to analyze the integral.
		The Marchenko-Pastur distribution is given by
		$$ f_{\lambda}(x) = \frac{1}{2\pi} \frac{\sqrt{(\lambda_{+} - x)(x - \lambda_{-})}}{\lambda x} $$
		where
		$ \lambda_{\pm} = (1 \pm \sqrt{\lambda})^2$
		and $0 < \lambda < 1$ is the ratio of the random matrix. The antiderivative of this function when weighted with $1/x$ happens to have a closed form
		
		\begin{align*}
			&\int \frac{\sqrt{(\lambda_{+} - x)(x - \lambda_{-})}}{x^2} dx =\\
			&=  \frac{1}{x} \sqrt{-\lambda^2 - (x-1)^2 + 2\lambda(1+x)} \\
			&+\arctan\left( \frac{1 + \lambda - x}{\sqrt{-\lambda^2 - (x-1)^2 + 2\lambda(1+x)}}\right) \\
			&- \frac{1+\lambda}{1-\lambda} \arctan\left( \frac{x + \lambda(2+x) - \lambda^2 -1}{(\lambda-1) \sqrt{\lambda^2 + 2\lambda(1+x) - (1+x)^2}}\right).
		\end{align*}
		Taking appropriate limits $x \rightarrow \pm \lambda_{\pm}$, we obtain
		$$ \int_{\lambda_{-}}^{\lambda_{+}} \frac{\sqrt{(\lambda_{+} - x)(x - \lambda_{-})}}{x^2} dx = \frac{2 \lambda \pi}{1 - \lambda}$$
		from which we deduce
		$$ \int_{0}^{\infty} \frac{ f_{N/D}(x)}{x} dx = \frac{1}{1 -  N D^{-1}} = \frac{D}{D-N}.$$
	\end{proof}

	\bibliographystyle{IEEEtran}
	\bibliography{references}

\end{document}